\documentclass{amsart}
\usepackage{amssymb}
\usepackage{xcolor}
\usepackage{needspace}

\theoremstyle{remark}
\newtheorem*{remark*}{Remark}

\usepackage{etoolbox}
\makeatletter
\patchcmd{\@settitle}
  {\uppercasenonmath\@title}
  {}{}{}
\patchcmd{\@setauthors}
  {\MakeUppercase{\authors}}
  {\authors}{}{}

\patchcmd{\section}
  {\scshape}
  {\upshape}{}{}

\patchcmd{\maketitle}
  {\uppercasenonmath\shorttitle}
  {}{}{}
\patchcmd{\maketitle}
  {\@nx\MakeUppercase{\the\toks@}}
  {\the\toks@}{}{}
\makeatother

\usepackage[
  colorlinks=true,
  urlcolor=black,
  citecolor=red,
  linkcolor=blue,
  pagebackref=true,
  linktocpage=true,
  pdfpagelabels=true,
  bookmarksnumbered=true,
  bookmarksopen=true
]{hyperref}

\hypersetup{pdftitle={Nondegeneracy of the normalized bubble for a fourth-order Liouville boundary problem},pdfauthor={Jiaoping Chen, Wenjing Chen and Shengbing Deng},pdfsubject={Weighted nondegeneracy for a fourth-order Liouville boundary problem},pdfkeywords={Biharmonic equation, Liouville boundary problem, nondegeneracy, conformal covariance, trace Adams inequality}}

\newtheorem{theorem}{Theorem}[section]

\newtheorem{lemma}[theorem]{Lemma}
\theoremstyle{remark}

\numberwithin{equation}{section}

\begin{document}

\title[Nondegeneracy for a Liouville boundary problem]{Nondegeneracy of the bubble for a fourth-order Liouville boundary problem}
\author[J. Chen]{Jiaoping Chen}
\address{School of Mathematics and Statistics, Southwest University, Chongqing, 400715, P.R. China}
\email{cjp123888@email.swu.edu.cn}
\author[W. Chen]{Wenjing Chen}
\address{School of Mathematics and Statistics, Southwest University, Chongqing, 400715, P.R. China}
\email{wjchen@swu.edu.cn}
\author[S. Deng]{Shengbing Deng$^{\dagger}$}
\address{School of Mathematics and Statistics, Southwest University, Chongqing, 400715, P.R. China}
\email{shbdeng@swu.edu.cn}
\thanks{$^{\dagger}$Corresponding author: Shengbing Deng.}
\thanks{The research was supported by the National Natural Science Foundation of China 12371121 and 12671143, and by the Fundamental Research Funds for the Central Universities SWU-KF26002.}
\makeatletter
\@namedef{subjclassname@2020}{2020 Mathematics Subject Classification}
\makeatother

\begin{abstract}
We establish the nondegeneracy of the bubble for
\[
\Delta^2u=0\quad \text{in }\mathbb R^4_+,\qquad
\partial_tu=0,\quad \partial_t\Delta u=e^{24\pi^2u}
\quad \text{on }\partial\mathbb R^4_+.
\]
Using Fourier analysis, stereographic projection and a
removable-singularity argument, we classify the weighted $L^2$
boundary kernel without assuming boundedness of the boundary
solution. This  determines the bounded half-space kernel.
Both kernels are generated by dilation and tangential translations.
We further obtain sharp estimates for the projected spherical inverse.
\end{abstract}

\maketitle

\section{Introduction}

We consider the fourth-order Liouville boundary problem
\begin{equation}\label{eq:nonlinear-problem}
\begin{cases}
\Delta^2u=0, & \text{in }\mathbb R^4_+,\\
\partial_tu=0, & \text{on }\partial\mathbb R^4_+,\\
\partial_t\Delta u=e^{24\pi^2u}, & \text{on }\partial\mathbb R^4_+,
\end{cases}
\end{equation}
where $\mathbb R^4_+=\{(x',t):x'\in\mathbb R^3,\ t>0\}$,
$\partial\mathbb R^4_+=\mathbb R^3\times\{0\}$, and
$-\partial_t$ denotes the outward normal derivative.

Problem \eqref{eq:nonlinear-problem} arises as the limiting boundary
problem in the blow-up analysis of the sharp trace Adams inequality
in dimension four \cite{ChenLuZhu2026}.
The trace Adams inequality is an exponential counterpart of sharp
Sobolev trace inequalities \cite{AcheChang2017,Escobar1988}
and is closely related to sharp Moser--Trudinger-type inequalities
\cite{Adams1988,Beckner1993}. For more related results, we refer the reader to  \cite{Case2020,Cianchi2008,NgoNguyenPhan2020} and the references therein.

Kernel classification is a basic ingredient in the linear analysis
of concentrating solutions; see, for example, \cite{DavilaDelPinoMusso2005}
for a problem with exponential Neumann data.
For conformally invariant equations, conformal transformations
provide a spectral approach to this classification.
D\'avila et al.~\cite{DavilaDelPinoSire2013} combined stereographic
projection with spherical harmonics to prove the nondegeneracy of
fractional Sobolev bubbles.
In the exponential setting, DelaTorre et al.~\cite{DelaTorreEtAl2026}
characterized the linearized kernel of a singular fractional
Liouville equation on $\mathbb R$ through a Steklov eigenvalue
problem on a bounded planar domain.
For higher-dimensional exponential equations,
Li et al.~\cite{LiLiWeiYang2026}
classified the bounded linearized kernel of the standard
$Q$-curvature bubble in $\mathbb R^n$, $n\ge3$.

In this paper, we classify the bounded kernel of the linearization
of \eqref{eq:nonlinear-problem} at the standard bubble.
To this end, we derive a boundary equation in $\mathbb R^3$
and prove that bounded half-space solutions are uniquely
determined by their traces through a biharmonic Poisson
representation \cite{ChenLuZhu2026,ChenZhang2025}.
We then classify the boundary kernel in a weighted $L^2$ space
by a direct spherical argument, without assuming boundedness
of the boundary solution.
In the three-dimensional case, this extends the bounded
kernel classification in
\cite[Corollary~1.1]{LiLiWeiYang2026}.
The Poisson representation then yields the nondegeneracy
in the class of bounded half-space solutions.

In order to state our results, we introduce some notations.   Set $a=(\pi/2)^{2/3}$ and consider the normalized
standard bubble
\[
\begin{aligned}
U(x',t)=-\frac{1}{8\pi^2}
\log\bigl((1+at)^2+a^2|x'|^2\bigr)+\frac{at}{4\pi^2\bigl((1+at)^2+a^2|x'|^2\bigr)},
\end{aligned}
\]
which solves \eqref{eq:nonlinear-problem} and satisfies
\[
U(0,0)=\sup_{\overline{\mathbb R^4_+}}U=0,
\qquad
\int_{\mathbb R^3}e^{24\pi^2U(x',0)}\,dx'=1.
\]
The dilation and tangential translation symmetries of
\eqref{eq:nonlinear-problem} generate the family
\begin{equation}\label{eq:bubble-family}
\begin{aligned}
\Psi_{\lambda,\xi}(x',t)
={}&\frac{1}{8\pi^2}\log\frac{2\lambda}{(2\pi^2)^{\frac{1}{3}}}
-\frac{1}{8\pi^2}
\log\bigl((1+\lambda t)^2+\lambda^2|x'-\xi|^2\bigr)\\
&+\frac{\lambda t}
{4\pi^2\bigl((1+\lambda t)^2+\lambda^2|x'-\xi|^2\bigr)}.
\end{aligned}
\end{equation}
Here $\lambda>0$, $\xi\in\mathbb R^3$ and $\Psi_{a,0}=U$.
Differentiating the bubble family at $(\lambda,\xi)=(a,0)$ gives
\begin{equation}\label{eq:jacobi-fields}
\Phi_0=\left.\partial_\lambda\Psi_{\lambda,0}\right|_{\lambda=a},
\qquad
\Phi_j=\left.\partial_{\xi_j}\Psi_{a,\xi}\right|_{\xi=0},
\quad j=1,2,3.
\end{equation}
To study the nondegeneracy, we linearize
\eqref{eq:nonlinear-problem} at $U$ and obtain
\begin{equation}\label{eq:linearized-problem}
\begin{cases}
\Delta^2\phi=0, & \text{in }\mathbb R^4_+,\\
\partial_t\phi=0, & \text{on }\partial\mathbb R^4_+,\\
\partial_t\Delta\phi=24\pi^2e^{24\pi^2U}\phi,
& \text{on }\partial\mathbb R^4_+.
\end{cases}
\end{equation}
Let $L_U\phi=\Delta^2\phi$, with domain consisting of functions
$\phi\in C^4(\mathbb R^4_+)\cap C^3(\overline{\mathbb R^4_+})$
satisfying the boundary conditions in \eqref{eq:linearized-problem}.
We write $\ker_{L^\infty}L_U=\ker L_U\cap L^\infty(\mathbb R^4_+)$.
{

For the corresponding boundary problem in $\mathbb R^3$, we set
\[
\rho_a(x)=e^{24\pi^2U(x,0)}=(1+a^2|x|^2)^{-3},
\qquad
\mathcal X_a=L^2(\mathbb R^3,\rho_a\,dx),
\]
and define the boundary operator
\[
\mathcal L_av=2(-\Delta)^{\frac{3}{2}}v-24\pi^2\rho_av.
\]
The operator $\mathcal L_a:\mathcal X_a\to\mathcal D'(\mathbb R^3)$
is well defined (see Lemma~\ref{lem:duality}), with the fractional
Laplacian defined by
\begin{equation}\label{eq:distributional-laplacian}
\big\langle(-\Delta)^{\frac{3}{2}}v,\eta\big\rangle
=\int_{\mathbb R^3}v(x)(-\Delta)^{\frac{3}{2}}\eta(x)\,dx,
\qquad \eta\in C_c^\infty(\mathbb R^3).
\end{equation}
Set
\[
Z_0(x)=\frac{1-a^2|x|^2}{1+a^2|x|^2},
\qquad
Z_j(x)=\frac{2ax_j}{1+a^2|x|^2},
\quad j=1,2,3.
\]}

\Needspace{10\baselineskip}
Our first result establishes the nondegeneracy of the bubble
and classifies the weighted boundary kernel.
\begin{theorem}\label{thm:kernel}
The kernel of $L_U$ in $C^4(\mathbb R^4_+)\cap C^3(\overline{\mathbb R^4_+})
\cap L^\infty(\mathbb R^4_+)$ is
\[
\ker_{L^\infty}L_U
=\operatorname{span}\{\Phi_0,\Phi_1,\Phi_2,\Phi_3\}.
\]
The kernel of $\mathcal L_a$ in $\mathcal X_a$ is
\begin{equation}\label{eq:weighted-kernel}
\ker_{\mathcal X_a}\mathcal L_a
=\operatorname{span}\{Z_0,Z_1,Z_2,Z_3\}.
\end{equation}
\end{theorem}

As a quantitative complement to Theorem~\ref{thm:kernel}, we derive
sharp inverse and coercivity estimates from the spherical spectrum.
Let $\mathcal H_\ell(\mathbb S^3)$ be the space of spherical harmonics
of degree $\ell$. Define the self-adjoint operator
$P_3:H^3(\mathbb S^3)\subset L^2(\mathbb S^3)
\to L^2(\mathbb S^3)$ by
\begin{equation}\label{eq:P3-spectrum}
P_3Y_\ell=\lambda_\ell Y_\ell,\qquad
\lambda_\ell=\ell(\ell+1)(\ell+2),\qquad
Y_\ell\in\mathcal H_\ell(\mathbb S^3).
\end{equation}
This is the standard conformally covariant operator of order three
on $\mathbb S^3$; see \cite{AcheChang2017,Branson1995}.
The spherical linearized operator is $P_3-6$, whose kernel is
$\mathcal H_1(\mathbb S^3)$.
For $F\in\mathcal D'(\mathbb S^3)$, define
$F_\ell\in\mathcal H_\ell(\mathbb S^3)$ by
\[
\int_{\mathbb S^3}F_\ell Y\,d\omega=\langle F,Y\rangle,
\qquad Y\in\mathcal H_\ell(\mathbb S^3),\quad \ell\ge0.
\]
The same convention applies to $W_\ell$.
Here $\mathcal D(\mathbb R^n)=C_c^\infty(\mathbb R^n)$ and $\langle\cdot,\cdot\rangle$ denotes the natural duality
pairing, which agrees with the usual $L^2$ inner product for real-valued $L^2$ functions.
We write $F\perp\mathcal H_1$ if $\langle F,Y\rangle=0$ for every
$Y\in\mathcal H_1$, or equivalently if $F_1=0$.
We use the Sobolev norms
\begin{equation}\label{eq:sobolev-norm}
\|W\|_{H^s(\mathbb S^3)}^2
=\sum_{\ell=0}^\infty(1+\ell)^{2s}
\|W_\ell\|_{L^2(\mathbb S^3)}^2.
\end{equation}
For every $s\in\mathbb R$, $P_3$ defines a bounded linear map
$P_3:H^s(\mathbb S^3)\rightarrow H^{s-3}(\mathbb S^3)$.
Set
\[
Q(W)=\big\langle(P_3-6)W,W\big\rangle,
\qquad W\in H^{\frac{3}{2}}(\mathbb S^3).
\]

The following theorem gives sharp inverse estimates for $P_3-6$
and a coercivity estimate for $Q$.
\begin{theorem}\label{thm:inverse}
Let $F\in H^{-3/2}(\mathbb S^3)$ satisfy
$F\perp\mathcal H_1(\mathbb S^3)$. Then there exists a unique
$W\in H^{3/2}(\mathbb S^3)\cap\mathcal H_1(\mathbb S^3)^\perp$
such that
\begin{equation}\label{eq:projected-equation}
(P_3-6)W=F.
\end{equation}
Moreover,
\begin{equation}\label{eq:inverse-estimates}
\|W\|_{H^{\frac{3}{2}}}\le\frac32\|F\|_{H^{-\frac{3}{2}}}.
\end{equation}
If $F\in L^2(\mathbb S^3)$, then $W\in H^3(\mathbb S^3)$ and
\begin{equation}\label{eq:strong-inverse}
\|W\|_{H^3}\le\frac32\|F\|_{L^2}.
\end{equation}
For $W\in H^{3/2}(\mathbb S^3)$ with $W\perp1$ and
$W\perp\mathcal H_1(\mathbb S^3)$,
\begin{equation}\label{eq:sharp-gaps}
Q(W)\ge18\|W\|_{L^2}^2,\qquad
Q(W)\ge\frac34\langle P_3W,W\rangle,\qquad
Q(W)\ge\frac23\|W\|_{H^{\frac{3}{2}}}^2.
\end{equation}
All constants in \eqref{eq:inverse-estimates}--\eqref{eq:sharp-gaps}
are sharp with respect to
the norms defined in \eqref{eq:sobolev-norm}.
For $W\ne0$, equality holds in any of the estimates in
\eqref{eq:sharp-gaps} if and only if
$W\in\mathcal H_2(\mathbb S^3)$.
The quadratic form $Q$ has Morse index one and nullity four on
$H^{3/2}(\mathbb S^3)$.
\end{theorem}

In Section~\ref{sec:reduction}, we establish the Poisson
representation and boundary reduction by Fourier analysis. Section~\ref{sec:kernel} combines stereographic projection
with spherical spectral analysis to prove
Theorems~\ref{thm:kernel} and~\ref{thm:inverse}.

\section[Boundary reduction]{{Boundary reduction}}
\label{sec:extension}\label{sec:reduction}

We use the Fourier convention
\[
\widehat f(\zeta)=\int_{\mathbb R^3}e^{-ix'\cdot\zeta}f(x')\,dx',
\qquad
f(x')=\frac{1}{(2\pi)^3}\int_{\mathbb R^3}e^{ix'\cdot\zeta}\widehat f(\zeta)\,d\zeta.
\]
For $t>0$ and $x'\in\mathbb R^3$, define
\begin{equation}\label{eq:poisson-kernels}
p_t(x')=\frac{t}{\pi^2(t^2+|x'|^2)^2},
\qquad
K_t(x')=\frac{4t^3}{\pi^2(t^2+|x'|^2)^3}.
\end{equation}
The kernels $p_t$ and $K_t$ are nonnegative and have unit mass.
Given $f\in L^\infty(\mathbb R^3)$, define $P_2(f)(x',t)=K_t*f(x')$. Let  $\mathcal S(\mathbb R^n)$ denote the Schwartz space.

{\begin{lemma}\label{lem:poisson}
For $f\in L^\infty(\mathbb R^3)$, $P_2(f)$ is smooth and
biharmonic in $\mathbb R^4_+$ and satisfies
\begin{equation}\label{eq:poisson-bound}
\|P_2(f)\|_{L^\infty(\mathbb R^4_+)}
\le \|f\|_{L^\infty(\mathbb R^3)}.
\end{equation}
If $f\in\mathcal S(\mathbb R^3)$, then
\begin{equation}\label{eq:bianjietiaojian}
P_2(f)|_{t=0}=f,\qquad
\partial_tP_2(f)|_{t=0}=0,
\end{equation}
and
\begin{equation}\label{eq:poisson-fourier}
\widehat{P_2(f)}(\zeta,t)
=(1+t|\zeta|)e^{-t|\zeta|}\widehat f(\zeta),
\qquad
\partial_t\Delta P_2(f)|_{t=0}
=2(-\Delta)^{\frac32}f.
\end{equation}
\end{lemma}}
\begin{proof}
By \eqref{eq:poisson-kernels}, $K_t=p_t-t\partial_tp_t$.
Since $\widehat p_t(\zeta)=e^{-t|\zeta|}$,
$\widehat K_t(\zeta)=(1+t|\zeta|)e^{-t|\zeta|}$.
For $k=|\zeta|$, we have
\[
(\partial_t^2-k^2)\bigl((1+kt)e^{-kt}\bigr)
=-2k^2e^{-kt},
\qquad
(\partial_t^2-k^2)^2\bigl((1+kt)e^{-kt}\bigr)=0.
\]
Hence $\Delta^2K_t=0$.
For $f\in L^\infty(\mathbb R^3)$, the \(L^1\)-bounds on the derivatives of \(K_t\), locally uniform for \(t>0\),
imply
$P_2(f)\in C^\infty(\mathbb R^4_+)$ and
\[
\Delta^2P_2(f)(x',t)
=\int_{\mathbb R^3}\Delta^2K_t(x'-y)f(y)\,dy=0.
\]
The nonnegativity and unit mass of $K_t$ imply
\eqref{eq:poisson-bound}.

For $f\in\mathcal S(\mathbb R^3)$,
$\widehat{P_2(f)}=\widehat K_t\,\widehat f$
gives the Fourier identity in \eqref{eq:poisson-fourier}.
Moreover,
\[
\bigl[(1+kt)e^{-kt}\bigr]_{t=0}=1,
\qquad
\left.\partial_t\bigl((1+kt)e^{-kt}\bigr)\right|_{t=0}=0,
\qquad
\partial_t\bigl(-2k^2e^{-kt}\bigr)=2k^3e^{-kt}.
\]
Since $\widehat f\in\mathcal S(\mathbb R^3)$,
Fourier inversion and dominated convergence as $t\to0^+$
yield the boundary identities
in \eqref{eq:bianjietiaojian}--\eqref{eq:poisson-fourier}.
\end{proof}


We next record the properties of the bubble family and its
symmetry fields.
\begin{lemma}\label{lem:bubbles}
The bubble family \eqref{eq:bubble-family} solves
\eqref{eq:nonlinear-problem}, and
the four fields \eqref{eq:jacobi-fields} are bounded solutions of
\eqref{eq:linearized-problem}. Moreover, 
\begin{equation}\label{eq:jacobi-traces}
\Phi_0(x',0)=\frac{Z_0(x')}{8\pi^2a},
\qquad
\Phi_j(x',0)=\frac{aZ_j(x')}{8\pi^2},\quad j=1,2,3.
\end{equation}
\end{lemma}
\begin{proof}
By \cite{ChenLuZhu2026}, the normalized bubble $U=\Psi_{a,0}$ solves
\eqref{eq:nonlinear-problem}. Since \eqref{eq:nonlinear-problem} is invariant under
translations and dilations
\[
u(x',t)\mapsto
u\bigl(\mu(x'-\xi),\mu t\bigr)+\frac{\log\mu}{8\pi^2},
\qquad \mu>0,
\]
the family $\Psi_{\lambda,\xi}$ also solves \eqref{eq:nonlinear-problem}.
Differentiating the equation and its boundary conditions with respect to
$\lambda$ and $\xi_j$, followed by evaluation at $(\lambda,\xi)=(a,0)$,
yields
\[
\Delta^2\Phi_i=0,\qquad
\left.\partial_t\Phi_i\right|_{t=0}=0,\qquad
\left.\partial_t\Delta\Phi_i\right|_{t=0}
=
24\pi^2e^{24\pi^2U(x',0)}\Phi_i,
\]
for $i=0,1,2,3$. Hence the fields in \eqref{eq:jacobi-fields} solve
\eqref{eq:linearized-problem}.
Set
\[
z=\lambda(x'-\xi),\qquad
\tau=\lambda t,\qquad
D=(1+\tau)^2+|z|^2,\qquad
E=\tau(1+\tau)+|z|^2.
\]
{Differentiation gives}
\begin{equation*}
8\pi^{2}\lambda \partial_{\lambda}\Psi_{\lambda,\xi}=1-\frac{2E}{D}+\frac{2\tau}{D}-\frac{4\tau E}{D^{2}},\qquad 4\pi^{2}\lambda^{-1}\partial_{\xi_{j}}\Psi_{\lambda,\xi}=\frac{z_{j}}{D}+\frac{2\tau z_{j}}{D^{2}}.
\end{equation*}
Since $
E\le D,$ $\tau\le\sqrt D,$ $|z|\le\sqrt D$ and $D\ge1$,
\[
\sup_{\mathbb R^4_+}|\partial_\lambda\Psi_{\lambda,\xi}|
\le C\lambda^{-1},
\qquad
\sup_{\mathbb R^4_+}|\partial_{\xi_j}\Psi_{\lambda,\xi}|
\le C\lambda.
\]
Thus the four fields in \eqref{eq:jacobi-fields} are bounded in $\mathbb{R}_+^4$. Evaluating the formulas at
$t=0$ and $(\lambda,\xi)=(a,0)$ yields \eqref{eq:jacobi-traces}.
\end{proof}


The boundary reduction uses the following decay estimate.
For $\eta\in C_c^\infty(\mathbb R^3)$,
\begin{equation}\label{eq:test-decay}
(-\Delta)^{\frac{3}{2}}\eta(x)=O((1+|x|)^{-6}).
\end{equation}
Indeed, write $(-\Delta)^{3/2}\eta=(-\Delta)^{1/2}(-\Delta\eta)$.
Choose $R\ge1$ with $\operatorname{supp}\eta\subset B_R$.
When $|x|\ge2R$, integration by parts in the kernel of $(-\Delta)^{1/2}$ gives
\[
|(-\Delta)^{\frac{3}{2}}\eta(x)|
\le C\int_{B_R}|\eta(y)|\,|x-y|^{-6}\,dy
\le C\|\eta\|_{L^1}|x|^{-6}.
\]

We now show that a bounded biharmonic function with zero normal
derivative is determined by its boundary trace.
\begin{lemma}\label{lem:bounded-extension}
Let $\phi\in C^4(\mathbb R^4_+)\cap C^3(\overline{\mathbb R^4_+})\cap L^\infty(\mathbb R^4_+)$ satisfy
\[
\Delta^2\phi=0\quad\text{in }\mathbb R^4_+,
\qquad \partial_t\phi=0\quad\text{on }\partial\mathbb R^4_+.
\]
Set $\phi_0(x')=\phi(x',0)$. Then
\begin{equation}\label{eq:bounded-extension}
\phi(x',t)=K_t*\phi_0(x'),\qquad t>0.
\end{equation}
Moreover, the boundary trace of $\partial_t\Delta\phi$ satisfies
\begin{equation}\label{eq:trace-flux}
\partial_t\Delta\phi|_{t=0}=2(-\Delta)^{\frac{3}{2}}\phi_0
\quad\text{in }\mathcal D'(\mathbb R^3).
\end{equation}
In particular, the trace of a bounded solution of \eqref{eq:linearized-problem} satisfies $\mathcal L_a\phi_0=0$.
\end{lemma}

\begin{proof}
Set $\widetilde\phi(x',t)=\phi(x',|t|)$ and
$g(x')=\left.\partial_t\Delta\phi(x',t)\right|_{t=0}$.
Since $\partial_t\phi|_{t=0}=0$, $\Delta\widetilde\phi$ is the even
extension of $\Delta\phi$, whose $t$-derivative has a jump $2g$ across
$t=0$. Hence
\begin{equation}\label{eq:even-distribution}
\Delta^2\widetilde\phi=2g\otimes\delta_0
\quad\text{in }\mathcal D'(\mathbb R^4),
\end{equation}
where $\delta_0$ is the Dirac mass at $t=0$ and
$\otimes$ denotes the tensor product of distributions.
{Choose $\theta\in C_c^\infty(\mathbb R)$ with $\theta(0)=1$.
For $\eta\in C_c^\infty(\mathbb R^3)$, \eqref{eq:even-distribution}
gives
\begin{equation*}
|\langle g,\eta\rangle|
=\frac12\big|\langle\widetilde\phi,\Delta^2(\eta\otimes\theta)\rangle\big|
\le C_\theta\|\phi\|_{L^\infty}
\sum_{|\alpha|\le4}\|\partial^\alpha\eta\|_{L^1}.
\end{equation*}
Thus $g\in\mathcal S'(\mathbb R^3)$, and
\eqref{eq:even-distribution} holds in $\mathcal S'(\mathbb R^4)$.
For $\chi\in C_c^\infty(\mathbb R^3\setminus\{0\})$ and $k=|\zeta|$,  the Fourier transform in $(x',t)$  yields}
\[
\chi(\zeta)\mathcal F_{x',t}\widetilde\phi(\zeta,\tau)
=\frac{2\chi(\zeta)\widehat g(\zeta)}{(k^2+\tau^2)^2}.
\]
{From
Fourier inversion in $\tau$ and the limit $t\to0^+$
in $\mathcal S'(\mathbb R^3)$, it follows that
\[
\chi\widehat\phi(\zeta,t)
=\frac{(1+kt)e^{-kt}}{2k^3}\chi\widehat g(\zeta)
=(1+kt)e^{-kt}\chi\widehat{\phi_0}(\zeta),
\qquad t>0.
\]}
{Let $u=K_t*\phi_0$. Then $u$ is bounded and biharmonic by
Lemma~\ref{lem:poisson}. Since $\chi\widehat K_t\in C_c^\infty$,
associativity of convolution gives
\[
\chi\widehat u(\zeta,t)
=(\chi\widehat K_t)\widehat{\phi_0}
=\bigl[\chi(1+kt)e^{-kt}\bigr]\widehat{\phi_0}(\zeta)
\quad\text{in }\mathcal S'.
\]
By the arbitrariness of $\chi$,}
\[
\operatorname{supp}
\widehat{\phi(\cdot,t)-u(\cdot,t)}
\subset\{0\}.
\]
{Thus $\phi(\cdot,t)-u(\cdot,t)$ is a polynomial in $x'$
and hence constant. 
By biharmonicity,
\[
\phi(x',t)-u(x',t)=c(t),\qquad
c^{(4)}=0,\qquad c\in L^\infty(0,\infty).
\]
Thus $c$ is constant. Since both $\phi(x',t)$ and
$u(x',t)=K_t*\phi_0(x')$ converge to $\phi_0(x')$ as $t\to0^+$,
we conclude that $c\equiv0$, proving \eqref{eq:bounded-extension}.}
{For $t>0$, the kernel $\partial_t\Delta K_t$ is even and integrable.
For $\eta\in C_c^\infty(\mathbb R^3)$, we differentiate
\eqref{eq:bounded-extension} and apply Fubini's theorem
and Lemma~\ref{lem:poisson} to obtain
\begin{equation*}
\langle\partial_t\Delta\phi(\cdot,t),\eta\rangle
=\langle\phi_0,(\partial_t\Delta K_t)*\eta\rangle=2\langle\phi_0,p_t*((-\Delta)^{\frac{3}{2}}\eta)\rangle.
\end{equation*}
By \eqref{eq:test-decay} and $(-\Delta)^{3/2}\eta\in L^1(\mathbb R^3)$, we have
\[
p_t*((-\Delta)^{\frac{3}{2}}\eta)\longrightarrow(-\Delta)^{\frac{3}{2}}\eta
\quad\text{in }L^1(\mathbb R^3)
\quad\text{as }t\to0^+.
\]
Since $\phi_0\in L^\infty$ and
$\phi\in C^3(\overline{\mathbb R^4_+})$, letting $t\to0^+$ yields
$\langle\partial_t\Delta\phi|_{t=0},\eta\rangle
=2\langle\phi_0,(-\Delta)^{3/2}\eta\rangle.$
This proves \eqref{eq:trace-flux}.} If $\phi$ solves \eqref{eq:linearized-problem}, then
\eqref{eq:trace-flux} implies $\mathcal L_a\phi_0=0$.
\end{proof}
\begin{remark*}
The boundedness assumption in Lemma~\ref{lem:bounded-extension}
is essential for the Poisson representation
\eqref{eq:bounded-extension}.
Indeed, $\phi(x',t)=t^2$ solves
\eqref{eq:linearized-problem} with
$\phi_0=0\in\mathcal X_a$, but
$\phi\ne K_t*\phi_0=0$.
\end{remark*}

\section[Kernel classification and spectral estimates]{{Kernel classification and spectral estimates}}
\label{sec:kernel}\label{sec:inverse}

To express the boundary equation on $\mathbb S^3$, we use
the stereographic projection
$\pi:\mathbb R^3\to\mathbb S^3\setminus\{N\}$,
where $N=(0,0,0,1)$, given by
\begin{equation}\label{eq:stereo}
\omega_j=\frac{2z_j}{1+|z|^2},\quad j=1,2,3,
\qquad \omega_4=\frac{|z|^2-1}{1+|z|^2},
\qquad \pi^*d\omega=\frac{8\,dz}{(1+|z|^2)^3}.
\end{equation}
For $v\in\mathcal X_a$, put $w(z)=v(z/a)$ and $W(\pi(z))=w(z)$.
By change of variables,
\begin{equation*}
\|W\|_{L^2(\mathbb S^3)}^2
=8a^3\int_{\mathbb R^3}\rho_a|v|^2\,dx
=2\pi^2\|v\|_{\mathcal X_a}^2<\infty.
\end{equation*}

\begin{lemma}\label{lem:duality}
{If $v\in\mathcal X_a$, then $v\in\mathcal S'(\mathbb R^3)$,
and the fractional Laplacian $(-\Delta)^{3/2}v$ is well defined
in $\mathcal D'(\mathbb R^3)$ by \eqref{eq:distributional-laplacian}.}
Let $W$ denote the spherical representative of $v$. Then
\begin{equation}\label{eq:dual-covariance}
\langle(P_3-6)W,\Theta\rangle_{\mathbb S^3}
=\frac12\langle\mathcal L_av,\Theta(\pi(a\,\cdot))\rangle_{\mathbb R^3},
\qquad \Theta\in C_c^\infty(\mathbb S^3\setminus\{N\}).
\end{equation}
\end{lemma}

\begin{proof}
Let $\eta\in\mathcal S(\mathbb R^3)$. Then 
\[
|\langle v,\eta\rangle|
\le\|v\|_{\mathcal X_a}
\left(\int_{\mathbb R^3}(1+a^2|x|^2)^3|\eta(x)|^2\,dx\right)^{\frac{1}{2}}.
\]
Thus $v\in\mathcal S'(\mathbb R^3)$.
If $\eta\in C_c^\infty(\mathbb R^3)$, \eqref{eq:test-decay} gives
\[
(-\Delta)^{3/2}\eta(x)=O((1+|x|)^{-6}).
\]
When $\operatorname{supp}\eta$ is contained in a fixed compact set,
the local Fourier bounds and the constant in
\eqref{eq:test-decay} are controlled by finitely many seminorms
of $\eta$.
Hence
\[
\int_{\mathbb R^3}(1+a^2|x|^2)^3|(-\Delta)^{\frac{3}{2}}\eta(x)|^2\,dx<\infty.
\]
By the Cauchy--Schwarz inequality,
the integral in \eqref{eq:distributional-laplacian}  is finite and depends continuously on \(\eta\).
Let $\Theta\in C_c^\infty(\mathbb S^3\setminus\{N\})$ and set
\[
\sigma(z)=\frac2{1+|z|^2},\qquad \eta(z)=\Theta(\pi(z)).
\]
{For $1<\gamma<3/2$, the spectral formula in \cite{Branson1995}
and the conformal covariance formula in
\cite{ChangGonzalez2011} give}
\begin{equation*}
P_{2\gamma} Y_\ell = \frac{\Gamma\left(\ell + \frac{3}{2} + \gamma\right)}{\Gamma\left(\ell + \frac{3}{2} - \gamma\right)} Y_\ell,\qquad (P_{2\gamma}\Theta)(\pi(z)) = \sigma(z)^{-\frac{3+2\gamma}{2}} (-\Delta)^\gamma \big(\sigma^{\frac{3-2\gamma}{2}}\eta\big)(z).
\end{equation*}
As $\gamma\to {3/2}^-$,
$\sigma^{(3-2\gamma)/2}\eta$ converges to $\eta$ in
$C_c^\infty(\mathbb R^3)$.
Uniform rapid decay of  $\widehat{\sigma^{(3-2\gamma)/2}\eta}$, together with
$|\zeta|^{2\gamma}\le1+|\zeta|^3$, yields
\[
(P_{2\gamma}\Theta)(\pi(z))
\longrightarrow
\sigma(z)^{-3}(-\Delta)^{\frac{3}{2}}\eta(z)
\quad\text{in }C^\infty_{\mathrm{loc}}(\mathbb R^3).
\]
For the spherical eigenvalues, we have
\[
0\le
\frac{\Gamma(\ell+3/2+\gamma)}
{\Gamma(\ell+3/2-\gamma)}
\le C(1+\ell)^3,
\qquad
\frac54\le\gamma<\frac32,
\]
and the ratio converges to $\lambda_\ell$. The decay of the harmonic
coefficients of $\Theta$ implies $
P_{2\gamma}\Theta\to P_3\Theta$ in $C^\infty(\mathbb S^3)$.
Thus $
(P_3\Theta)(\pi(z))
=
\sigma(z)^{-3}(-\Delta)^{3/2}\eta(z).$
Using the self-adjointness of $P_3$ and \eqref{eq:stereo},
we obtain
\[
\langle P_3W,\Theta\rangle_{\mathbb S^3}
=
\int_{\mathbb R^3}w(-\Delta)^{\frac{3}{2}}\eta\,dz,
\qquad
\langle W,\Theta\rangle_{\mathbb S^3}
=
8\int_{\mathbb R^3}
\frac{w\eta}{(1+|z|^2)^3}\,dz.
\]
Set $
\eta_a(x)=\eta(ax)=\Theta(\pi(ax)).$
Changing the  variable $z=ax$ gives 
\begin{equation*}
\langle(P_3-6)W,\Theta\rangle_{\mathbb S^3}
=
\int_{\mathbb R^3}v(-\Delta)^{\frac{3}{2}}\eta_a\,dx
-48a^3\int_{\mathbb R^3}\rho_av\eta_a\,dx\\
=\frac12\langle\mathcal L_av,\eta_a\rangle_{\mathbb R^3}.
\end{equation*}
\end{proof}

{The next lemma shows that an $L^2(\mathbb S^3)$ solution
of $(P_3-6)W=0$ away from $N$ satisfies the equation
on the whole sphere.}

\begin{lemma}\label{lem:removable}
Let $W\in L^2(\mathbb S^3)$ satisfy
 $
(P_3-6)W=0$ in $\mathcal D'(\mathbb S^3\setminus\{N\}).$
Then
\[
(P_3-6)W=0
\quad\text{in }\mathcal D'(\mathbb S^3).
\]
\end{lemma}
\begin{proof}
Set $T=(P_3-6)W$. By \eqref{eq:P3-spectrum},
\[
\|T\|_{H^{-3}}^2
=\sum_{\ell\ge0}(1+\ell)^{-6}|\lambda_\ell-6|^2
\|W_\ell\|_{L^2}^2
\le C\|W\|_{L^2}^2.
\]
Thus $T\in H^{-3}(\mathbb S^3)$. Moreover,
$\operatorname{supp}T\subset\{N\}$.
In local coordinates centered at \(N\), the structure theorem for point-supported distributions \cite[Theorem~2.3.4]{Hormander2003} shows that \(T\) is a finite linear combination of derivatives of \(\delta_0\).
Since
$\widehat{\partial^\alpha\delta_0}=(i\zeta)^\alpha$, we obtain
\begin{equation*}
\partial^\alpha\delta_0\in H^{-s}(\mathbb R^3)
\quad\Longleftrightarrow\quad
s>\frac32+|\alpha|.
\end{equation*}
{By $T\in H^{-3}$, we obtain}
\[
T=c_0\delta_N+\sum_{j=1}^3c_j\partial_j\delta_N.
\]
Choose geodesic coordinates at $N$ with tangent basis given by the first
three coordinate vectors of $\mathbb R^4$. Then
\[
\omega_4(N)=1,\qquad \partial_j\omega_4(N)=0,\qquad
\omega_i(N)=0,\qquad \partial_j\omega_i(N)=\delta_{ij},
\quad 1\le i,j\le3.
\]
Since $(P_3-6)\omega_i=0$ for $i=1,2,3,4$, self-adjointness yields
$$c_0 = \langle T, \omega_4 \rangle = \langle W, (P_3 - 6)\omega_4 \rangle = 0$$ and
\begin{equation*}
-c_i = \langle T, \omega_i \rangle = \langle W, (P_3 - 6)\omega_i \rangle = 0,\quad i=1,2,3.
\end{equation*}
Thus $T=0$.
\end{proof}

We now classify the weighted boundary kernel and determine
the corresponding bounded half-space solutions.

\begin{proof}[Proof of Theorem~\ref{thm:kernel}]
Let $v\in\mathcal X_a$ solve $\mathcal L_av=0$.
Lemmas~\ref{lem:duality} and~\ref{lem:removable} yield
\[
(P_3-6)W=0\quad\text{in }\mathcal D'(\mathbb S^3).
\]
By the spectral decomposition \eqref{eq:P3-spectrum}, for each $\ell\ge0$, $(\lambda_\ell-6)W_\ell=0$.
Since $\lambda_\ell=6$ if and only if $\ell=1$, we obtain
 $W\in\mathcal H_1(\mathbb S^3)=\operatorname{span}\{\omega_1,\omega_2,\omega_3,\omega_4\}$.
Consequently, \eqref{eq:stereo} and $v(x)=W(\pi(ax))$ give
\[
v\in\operatorname{span}\{Z_0,Z_1,Z_2,Z_3\}.
\]
Conversely, $Z_0,Z_1,Z_2,Z_3\in\mathcal X_a$ have spherical
representatives $-\omega_4,\omega_1,\omega_2,\omega_3$, respectively.
Since $(P_3-6)\omega_i=0$, $i=1,2,3,4$,
\eqref{eq:dual-covariance} implies $\mathcal L_aZ_j=0$
for $j=0,1,2,3$. This proves \eqref{eq:weighted-kernel}.
Now let $\phi\in\ker_{L^\infty}L_U$. Lemma~\ref{lem:bounded-extension}
gives $\phi_0\in\mathcal X_a$ and $\mathcal L_a\phi_0=0$.
By \eqref{eq:weighted-kernel} and \eqref{eq:jacobi-traces},
there exist $d_0,d_1,d_2,d_3$ such that
\[
h=\phi-\sum_{j=0}^3d_j\Phi_j
\]
has zero boundary trace. Lemma~\ref{lem:bubbles} ensures that
$h$ is bounded and satisfies
\[
\Delta^2h=0 \quad\text{in }\mathbb R^4_+,
\qquad
h=\partial_th=0 \quad\text{on }\partial\mathbb R^4_+.
\]
From Lemma~\ref{lem:bounded-extension}, we obtain
\[
h(x',t)=K_t*h(\cdot,0)(x')=0.
\]
Hence
$\phi\in\operatorname{span}\{\Phi_0,\Phi_1,\Phi_2,\Phi_3\}$.
The reverse inclusion follows from Lemma~\ref{lem:bubbles}.
\end{proof}

We now prove the inverse and coercivity estimates.

\begin{proof}[Proof of Theorem~\ref{thm:inverse}]
If \eqref{eq:projected-equation} holds, self-adjointness gives
\[
\langle F,Y\rangle=\langle W,(P_3-6)Y\rangle=0,
\qquad Y\in\mathcal H_1.
\]
Assume
$F\perp\mathcal H_1$, so that $F_1=0$, and define
\begin{equation}\label{eq:spectral-inverse}
W=-\frac16F_0+\sum_{\ell=2}^\infty
\frac{F_\ell}{\lambda_\ell-6}.
\end{equation}
For $\ell\ge2$, we have
\begin{equation}\label{eq:spectral-gap}
3(\lambda_\ell-6)-2(1+\ell)^3
=(\ell-2)(\ell^2+5\ell+10)\ge0,
\end{equation}
with equality only at $\ell=2$. Hence, for $s\in\{-3/2,0\}$ and
$F\in H^s$, we  obtain
\begin{equation*}
\frac1{36}\|F_0\|_{L^2}^2+
\sum_{\ell=2}^\infty
\frac{(1+\ell)^{2s+6}}{(\lambda_\ell-6)^2}
\|F_\ell\|_{L^2}^2
\le\frac94\|F\|_{H^s}^2.
\end{equation*}
This proves convergence of \eqref{eq:spectral-inverse} in
$H^{s+3}$ and yields \eqref{eq:inverse-estimates}
and \eqref{eq:strong-inverse}.
The function $W$ defined by \eqref{eq:spectral-inverse} satisfies \eqref{eq:projected-equation} and
$W\perp\mathcal H_1$.
Since $\ker(P_3-6)=\mathcal H_1$, $W$ is the unique solution
{of \eqref{eq:projected-equation}}
in $H^{3/2}\cap\mathcal H_1^\perp$.
For $0\ne F\in\mathcal H_2$, equality holds in
\eqref{eq:inverse-estimates} and \eqref{eq:strong-inverse},
so the constant $3/2$ is sharp in each estimate.
If $W\perp1$ and $W\perp\mathcal H_1$, then $W_0=W_1=0$ and
\[
Q(W)=\sum_{\ell=2}^\infty
(\lambda_\ell-6)\|W_\ell\|_{L^2}^2.
\]
For $\ell\ge2$, the bound $\lambda_\ell\ge24$ and
\eqref{eq:spectral-gap} imply
\[
\lambda_\ell-6\ge18,\qquad
\lambda_\ell-6\ge\frac34\lambda_\ell,\qquad
\lambda_\ell-6\ge\frac23(1+\ell)^3,
\]
with equality only at $\ell=2$.
{This proves \eqref{eq:sharp-gaps}.}
Finally, for  $W\in H^{3/2}(\mathbb S^3)$,
\[
Q(W)=-6\|W_0\|_{L^2}^2+
\sum_{\ell=2}^\infty
(\lambda_\ell-6)\|W_\ell\|_{L^2}^2.
\]
The negative spectral subspace is $\mathcal H_0$, of dimension one,
and the kernel of the associated bilinear form is $\mathcal H_1$,
of dimension four. Thus $Q$ has Morse index one and nullity four.
\end{proof}

\noindent\textbf{Data availability.}
No data were used for the research described in this article.

\vspace{.1cm}

\noindent\textbf{Conflict of interest.}
The authors declare that they have no conflict of interest.

\end{document}